\documentstyle[12pt]{article}

\newtheorem{Thm}{Theorem}[section]

\newtheorem{Cor}[Thm]{Corollary}
\newtheorem{Lem}[Thm]{Lemma}

\newtheorem{Def}[Thm]{Definition}
\newtheorem{Rmk}[Thm]{Remark}

\newcommand{\im}{\mbox{$\Rightarrow$}}

\newcommand{\e}{\varepsilon}

\newcommand{\subl}{{_{\mbox {\scriptsize $\ell$}}}}

\newcommand{\sub}[1]{{_{\mbox {\scriptsize $#1$}}}}
\newcommand{\spn}{{\rm span}\,}

\newcommand{\ext}{{\rm ext}\,}
\newcommand{\next}{{\rm next}\,}
\newcommand{\wstarcl}{{\rm w}^*{\rm -cl}\,}
\newcommand{\wstarap}{{\rm w}^*{\rm -ap}\,}

\def \norm{| \! | \! |}
\def \N{\rm {\bf N}}
\def \R{\rm {\bf R}}

\begin{document}

\title{A subsequence characterization of sequences spanning
isomorphically polyhedral Banach spaces}

\author{G. Androulakis}

\date{October 13, 1996}
\maketitle

\noindent
{\bf Abstract:}
Let $(x_n)$ be a sequence in a Banach space $X$ which does not
converge in norm, and let $E$ be an isomorphically precisely norming
set for $X$ such that 
\[ \sum_n |x^*(x_{n+1}-x_n)|< \infty, \; \forall x^* \in E.  
\qquad (*) \]
Then there exists a subsequence of $(x_n)$ which spans an
isomorphically polyhedral Banach space. It follows immediately from
results of V. Fonf that the converse is also true: If a separable
Banach space $Y$ is a separable isomorphically polyhedral then there
exists a non
norm convergent sequence $(x_n)$ which spans $Y$ and there exists an
isomorphically precisely norming set $E$ for $Y$ such that $(*)$ is
satisfied. As an application of this subsequence characterization of
sequences spanning isomorphically polyhedral Banach spaces we obtain a
strengthening of a result of J. Elton, and an Orlicz-Pettis type result.

\bigskip
\noindent
{\bf Acknowledgments:}
I would like to thank Professor H. Rosenthal for suggesting this
project to me and for his help. Also, I would like to thank Professors
N. J. Kalton and G. Godefroy for their comments.

\noindent
\section{Introduction} \label{S:intro}
In 1958 C. Bessaga and A. Pelczynski proved the following 

\begin{Thm} [\cite{BP}] \label{T:bp}
If $(x_n)$ is a non-weakly convergent sequence
in a Banach space $X$ such that 
\begin{equation} \label{E:bp}
\sup_{x^* \in Ba(X^*)} \sum_n | x^*(x_{n+1}-x_n) |< \infty 
\end{equation}
then there exists a subsequence of $(x_n)$ which is equivalent to the
summing basis $(s_n)$ of $c_0$.
\end{Thm}

\noindent
Recall that the summing basis $(s_n)$ of $c_0$ is defined by $s_n= e_1
+ \ldots + e_n, \; \forall n \in \N$, where $(e_n)$ denotes the unit
vector basis of $c_0$. In 1981 J. Elton was able to eliminate the
assumption ``non-weakly convergent'' and relax the condition
(\ref{E:bp}) and still obtain that $c_0$ embeds in the closed linear
span $[x_n]$ of $(x_n)$. The result of J. Elton can be stated as
follows:

\begin{Thm} [\cite{E2}] \label{T:elton}
If $(x_n)$ is a semi-normalized basic sequence in a Banach space $X$
such that 
\begin{equation} \label{E:elton}
\sum |x^* (x_n)| < \infty, \; \forall x^* \in \ext Ba(X^*)
\end{equation}
(where $\ext Ba(X^*)$ denotes the set of the extreme points of the
dual ball) then $c_0$ embeds in $[x_n]$. 
\end{Thm}

\noindent
In order to prove this result, J. Elton first showed that
there exists a polyhedral Banach space which embeds in $[x_n]$  (for
the definition, examples and properties of the polyhedral Banach
spaces see the next section). Then the result of Theorem \ref{T:elton}
follows from the following theorem of V. Fonf:

\begin{Thm} [\cite{F3}]
Every polyhedral Banach space $X$ contains an isomorph of $c_0$, and
if in addition $X$ is separable, then $X^*$ is separable.
\end{Thm}

\noindent
We prove a stronger result than Theorem \ref{T:elton} by eliminating
the condition of having a basic sequence, by replacing the
set of the extreme points in condition (\ref{E:elton}) by any
isomorphically precisely norming set, and finally by obtaining the precise way
that a polyhedral Banach space embeds in $[x_n]$. 
Our main result can be
stated as follows:

\begin{Thm} \label{T:main}
If $(x_n)$ is a sequence in a Banach space $X$ which does not converge
in norm, 
and $E$ is an isomorphically precisely norming set for $X$ such that 
\begin{equation} \label{E:main}
\sum_n | x^*(x_{n+1}-x_n)| < \infty, \; \forall x^* \in E 
\end{equation}
then there exists a subsequence of $(x_n)$ which spans an
isomorphically polyhedral Banach space.

\noindent
Conversely, if $Y$ is a separable isomorphically polyhedral Banach
space then there exists a non norm convergent sequence $(x_n)$ in $Y$
and an isomorphically precisely norming set $E$ for $Y$ such that
$[x_n]=Y$ and (\ref{E:main}) holds.
\end{Thm}

\noindent
Recall that a set $E
\subset X^*$ is called {\it isomorphically precisely norming for $(X,
\| \cdot \|)$}, (the terminology is due to H. Rosenthal \cite{R}), if
there  exists $C \geq 1$ such that 
\begin{itemize}
\item[(a)] $E \subseteq C \cdot Ba(X^*)$,
\item[(b)] $\frac{1}{C} \| x \| \leq \sup_{e \in E} |e(x)|, \; \forall
         x \in X$, and
\item[(c)] $\forall x \in X \; \exists e_0 \in E \; |e_0(x)|= \sup_{e
\in E} |e(x) |$.
\end{itemize}
If $E$ satisfies (a), (b), and (c) for $C=1$ then $E$ is called {\it
precisely norming for $(X, \| \cdot \|)$}. 

\begin{Def}
The set of the norm achieving extreme points of the dual ball of a
Banach space $X$ is defined as follows:
\[ \next Ba(X^*)= \{ x^* \in \ext
Ba(X^*) : \exists x \in Ba(X) \; |x^*(x)|=1 \}. \] 
\end{Def}

\noindent
The set $\next Ba(X^*)$ is an example of a precisely norming set for
$X$.

\bigskip
\noindent
Theorem \ref{T:main} is a strengthening of the following remark which
can be easily derived from a result of V. Fonf \cite{F4}. 

\begin{Rmk} \label{R:fonf}
Under the same hypotheses of Theorem \ref{T:main} there exist a
sequence $(\e _n) \in \{ \pm 1 \} ^{\N}$ and an increasing sequence of
positive integers $(\ell_k)$ such that $[( \sum_{i=1}^{\ell_k} \e_i
(x_i - x_{i-1}))_k]$ is an i.p. space.
\end{Rmk}

\noindent
We sketch the proof of Remark \ref{R:fonf} at the end of Section \ref{S:proof}.

\bigskip
\noindent
The last section is devoted to applications of Theorem \ref{T:main}. 
One application is given in $C(K)$ spaces.
If $K$ is a compact metric space then $DSC(K)$ denotes the class of
bounded differences of semi-continuous functions on $K$ (the
definition appears in section \ref{S:applications}). An immediate
corollary of Theorem \ref{T:main} is the following:

\begin{Thm} \label{T:dsc}
Let $f \in DSC(K) \backslash C(K)$ be given, where  $K$ is a compact
metric  space. Then $f$ strictly governs the class of (separable)
polyhedral  Banach spaces.
\end{Thm}

\noindent
This theorem was the main motivating result for this research.
The definitions of the terms ``strictly governs'' and ``governs''
appear in section \ref{S:applications}. This generalizes the following theorem
of J. Elton which was also proved by R.Haydon, E. Odell and H. Rosenthal:

\begin{Thm} [\cite{E2}, \cite{HOR}] \label{T:hordsc}
Let $f \in DSC(K) \backslash C(K)$ be given, where $K$ is a compact
metric  space. Then $f$ governs $\{ c_0 \}$.
\end{Thm}

\noindent
Another application is the following Orlicz-Pettis type result:

\begin{Thm} \label{T:orliczpettis}
Let $(y_n)$ be a sequence in a Banach space $X$ and let $E$ be an
isomorphically precisely norming set for $X$. If $c_0$ does not embed
isomorphically in the closed linear span $[y_n]$ of $(y_n)$ and 
\[ \sum_n | x^*(y_n) | < \infty , \; \forall x^* \in E, \]
then $\sum_n y_n$ converges unconditionally.
\end{Thm}
 
\section{Isomorphically polyhedral Banach spaces} \label{S:ip}

\noindent
Polyhedral Banach spaces were introduced by V. Klee \cite{K}. An
infinite dimensional Banach space is called {\it polyhedral} if the ball of
any of its finite dimensional subspaces is a polyhedron, i.e. it has
finitely many extreme points. $c_0$ is an example of a polyhedral
Banach space. A Banach space will be called {\it  isomorphically
polyhedral} (i.p. in short) if it is polyhedral under some equivalent
norm. We are interested in isomorphic theory and therefore in i.p.
Banach spaces. Examples of i.p. Banach spaces are: $c$ (the space of
convergent sequences), the $\ell_1$ preduals \cite{F4}, the spaces
$C(\alpha)$ for any ordinal $\alpha$ \cite{F2}, $c_0$-sum of
separable i.p. spaces (easy to prove using Theorem
\ref{T:charactpolyh}), finite dimensional extensions of i.p. spaces
(easy to prove), the Orlicz sequence space $h_M$ where $M$ is a
non-degenerate Orlicz function satisfying $\lim_{t \rightarrow
0}M(Kt)/M(t)= \infty$ for some $K >1$ \cite{L}. The following
characterization of the separable i.p. Banach spaces was proved by V.
Fonf: (note that if $(X, \| \cdot \|)$ is a Banach space, $| \cdot |$
is an equivalent norm and $C \geq 1$ then we say that these norms are
{\it $C$ equivalent} if $C^{-1} \| x \| \leq | x | \leq C \| x
\|$ for all $x \in X$)

\begin{Thm} [\cite{F3}, \cite{F4}, \cite{F5}] \label{T:charactpolyh}
Let $(X,\| \cdot \|)$ be a separable Banach space. TFAE
\begin{itemize}
\item[(1)] For every $\e >0$ there exists a $1+ \e$ equivalent norm $|
\cdot |$ on $X$ such that $(X, | \cdot |)$ is polyhedral.
\item[(2)] For every $\e >0$ there exists a $1 + \e$ equivalent norm
$\norm \cdot \norm$ 
on $X$ such that the set $\next Ba(X, \norm \cdot \norm)^*$ is countable.
\end{itemize}
\end{Thm}

\noindent
The following theorem which was proved recently by R. Deville, V. Fonf
and P. H\'{a}jek, gives an interesting property of the i.p. Banach
spaces.

\begin{Thm} [\cite{DFH}] \label{T:dfh}
Let $(X, \| \cdot \|)$ be a separable i.p. Banach space, and $\e >0$.
Then there exists a $1 + \e$ equivalent norm $| \cdot |$ on $X$ such
that $(X, | \cdot |)$ is polyhedral.
\end{Thm}

\noindent
The following lemma, which gives a sufficient condition for a Banach
space $X$ to be an i.p. space, uses ideas similar to the ones that are
found in \cite{F1}. It is stated and proved in a slightly more general
setting that we shall need. Before we state it, we need some
terminology: Let $X$ be a Banach space and $K$ be a weak$^*$ compact
subset of $X^*$. The Banach space $C(K)$ will always be equipped with
the supremum norm which will be denoted by $\| \cdot \|_\infty$ or $\|
\cdot \|_{C(K)}$. We denote by $X \mid K$ the subspace of $C(K)$ which
is obtained as the image of $X$ via the composition of the maps:
\[ X \stackrel{\phi}{\longrightarrow} X^{**}
\stackrel{r}{\longrightarrow} C(K)\]
where $\phi$ is the canonical embedding of $X$ in $X^{**}$ and $r$ is
the restriction map defined by $r(x^{**})= x^{**} \mid K$ for every
$x^{**} \in X^{**}$ (where the domain of $x^{**} \mid K$ is $K$, and 
\[ (x^{**} \mid K)(k)= x^{**}(k), \mbox{ for all } k \in K).\] 
Note that $X \mid K$ is a (possibly non-complete) normed linear space.
For $x \in X$ we write $x \mid K= r(\phi (x))$.
We adopt the convention that the maximum or the supremum of an empty
family of numbers is $- \infty$ and the minimum or the infimum of an
empty family of numbers is $+ \infty$.

\begin{Lem} \label{L:ip}
Let $(X, \| \cdot \|)$ be a Banach space and $K \subset X^*$ be a
precisely norming set. Let $I, J$ be disjoint subsets of the integers.
Let $(K_n)_{n \in I \cup J}$ be weak$^*$ compact subsets of $K$ with
$K= \cup_{n \in I \cup J}K_n$ (one of the index sets $I, J$ may also
be empty). Let $(C_i)_{i \in I} \subset(1,2)$ and $(\e_j)_{j \in J}
\subset (0,1)$ be two sequences of numbers such that $C_i \searrow 1$
(if $I$ is infinite) and $\e_j \searrow 0$ (if $J$ is infinite). We
assume the following:
\begin{itemize}
\item For every $i \in I$ there exists an equivalent norm $\| \cdot
\|_i$ on $X \mid K_i$ such that 
\begin{eqnarray*}
& & \| y \|_\infty \leq \| y \|_i \leq C_i \| y \|_\infty, \; \forall
 y  \in X \mid K_i \mbox{ and} \\
& & (X \mid K_i, \| \cdot \|_i) \mbox{ has a countable precisely
 norming set}\\
& & \mbox{(we call it $\tilde{K}_i$ for future reference during the
 proof).}
\end{eqnarray*}
\item For every $j \in J$ there exists a finite set $A_j$ of $K_j$
such that for every $x \in X$,
\[ \mbox{if \ } \| x \|= \| x \mid  K_j \|_\infty \mbox{ then } \| x \|
\leq \frac{1}{1- \e_j} \max_{a \in A_j} | a(x) |. \]
\end{itemize}
Then $X$ is an i.p. space.
\end{Lem}

\noindent
{\bf Proof }  If $J \not = \emptyset$ then choose a sequence
$(\delta_j)_{j \in J}$ of positive numbers satisfying the following
properties:
\begin{itemize}
\item $1 +\delta_j > (\sup_{i \in I, i > j}C_i) \frac{1}{1-\e_j}$ for
all $j \in J$ with $\{ i \in I: i >j \} \not = \emptyset$,
\item $1 + \delta_j > \frac{1}{1-\e_j}$, (if $\{ i \in I: i > j \} =
\emptyset$), and
\item $\lim_{J \ni j \rightarrow \infty} \delta_j = 0$ (if $J$ is
infinite).
\end{itemize}
We can choose $(\delta_j)_{j \in J}$ with the above properties since
$C_i \searrow 1$ (if $I$ is infinite) and $\e_j \searrow 0$ (if $J$ is
infinite). Let $(x_{i_\subl , \ell}^*)_{\ell \in \N}$ be an
enumeration of the set $\cup_{i \in I} \tilde{K}_i$
such that 
\begin{eqnarray*}
& & i_\subl \in I \mbox{  and}\\ 
& &  x_{i_\subl , \ell}^* \in \tilde{K}_i
\end{eqnarray*}
for all $\ell \in \N$. For every $i \in I$ we consider the operator 
\[ T_i : (X, \| \cdot \|) \longrightarrow (C(K_i), \| \cdot \|_\infty)\]
defined by $T_i = r_i \circ \phi$ where $\phi : X \longrightarrow
X^{**}$ is the canonical embedding and $r_i :X^{**} \longrightarrow
C(K_i)$ is the restriction operator ($r_i(x^{**})=x^{**} \mid K_i$ for
all $x^{**} \in X^{**}$). Then, we define an equivalent norm $\norm
\cdot \norm$ on $X$ by 
\[ \norm x \norm = \sup_{\ell \in \N} (1 + \frac{1}{\ell}) |
(T_{i_\subl}^* x_{i_\subl, \ell}^*)x| \vee \sup_{j \in J}(1 +
\delta_j) \max_{a \in A_j} | a(x) |. \]
We shall show that $\norm \cdot \norm$ is an equivalent norm on $X$
such that $\next Ba(X, \norm \cdot \norm)^*$ is countable. This will
finish the proof of the lemma. 

\bigskip
\noindent
We first show that $\norm \cdot \norm$ is an equivalent norm on $X$. 

\bigskip
\noindent
Fix $x \in X$ with $\| x \|=1$. We have that $\| T_i \| \leq 1$ hence
$\| T_i^* \| \leq 1$ for all $i \in I$. Also we have that $\|
x_{i_\subl, \ell}^* \|_{C(K_{i_\subl})} \leq C_{i_\subl} \|
x_{i_\subl, \ell}^* \|_{i_\subl}= C_{i_\subl} <2$ for every $\ell \in
\N$. Thus we obtain 
\[ \norm x \norm \leq 4 \vee (1 + \max_{j \in J}
\delta_j). \] 
On the other hand there exists $k \in K$ such that
$|k(x)|=1$. Since $K= \cup_{n \in I \cup J}K_n$, there exists $n \in I
\cup J$ such that $k \in K_n$, and therefore $1= \| x \|= \| x \mid
K_n \|_\infty$. If $n \in I$, then since the set $\tilde{K}_i$
is precisely norming for $(X \mid K_n, \| \cdot
\|_n)$, there exists $\ell \in \N$ such that $i_\subl =n$ and 
\[ \| x \mid K_n \|_\infty \leq \| x \mid K_n \|_n = | x_{i_\subl ,
\ell} ^*(x \mid K_n)| = |( T_n^* x_{i_\subl, \ell}^*)(x) | \leq
\norm x \norm. \]
If $n \in J$ then 
\[1 = \| x \| \leq \frac{1}{1-\e_n} \max _{a \in A_n} | a (x) | \leq
(\max_{j \in J} \frac{1}{1-\e_j}) \norm x \norm. \]
Thus
\[ \norm x \norm \geq \min_{j \in J} (1-\e_j) \wedge 1. \]
Thus, $\norm \cdot \norm$ is an equivalent norm on $X$.

\bigskip
\noindent
For $x \in X$ we define 
\[ | x | = \sup_{\ell \in \N} | ( T_{i_\subl}^*x_{i_\subl,
\ell}^*)(x)| \vee \sup_{j \in J} \max _{a \in A_j} |a(x) |. \]
Then $| \cdot |$ is an equivalent norm on $X$, since it is evidently
equivalent to $\norm \cdot \norm$. Note that since
for every $i \in I$, $\tilde{K}_i$ is a
precisely norming set for $(X \mid K_i, \| \cdot \|_i)$, we have that 
\[ \sup_{\ell \in \N} | T_{i_\subl}^* x_{i_\subl, \ell}^*(x) |=
\sup_i \| x \mid K_i \|_i, \; \forall x \in X.\]

\bigskip
\noindent
Obviously, we have that $| x | \leq \norm x \norm$ for all $x \in X$.
We now show that $| x | < \norm x \norm$ for all $x \in X \backslash
\{ 0 \}$.

\bigskip
\noindent
Fix $x \in X$. Since $K$ is a precisely norming set for $(X, \| \cdot
\|)$, there exists $k \in K$ such that $\| x \| = | k(x) |$. Since $K
= \cup _{n \in I \cup J} K_n$, there exists $n \in I \cup J$ such that
$k \in K_n$; thus $\| x \| = \| x \mid K_n \|_\infty $. We separate
two cases:

\bigskip
\noindent
\underline{Case 1:} Assume that $n \in I$.

\bigskip
\noindent
There exists $\ell \in \N$ such that $i_\subl = n$ and $\| x \mid K_n
\|_n = |x_{i_\subl, \ell}^*(x \mid K_n)| = |(T_{i_\subl}^*x_{i_\subl
, \ell}^*)(x)|$. Then we have that 
\begin{eqnarray*}
\sup_{j \in J} \max_{a \in A_j} |a(x)| & \leq & \sup_{j \in J} \| x
\mid K_j \|_\infty \leq \| x \mid K_n \|_\infty \leq \| x \mid K_n
\|_n \\
& = & | (T_{i_\subl}^* x_{i_\subl , \ell}^*)(x) | < (1 +
\frac{1}{\ell}) | (T_{i_\subl}^* x_{i_\subl, \ell}^*)(x) | \leq \norm
x \norm. 
\end{eqnarray*}
Now, choose $i' \in I$ such that $\sup_{i \in I, i > i'}C_i < 1 +
\frac{1}{\ell}$. Hence 
\begin{eqnarray*}
\sup_{i \in I, i > i'} \| x \mid K_i \|_i & \leq & \sup_{i \in I, i
>i'}C_i \| x \mid K_i \|_\infty \\
& \leq & \sup _{i \in I, i > i'}C_i \| x \mid K_n \|_\infty \\
& \leq & \sup_{i \in I, i >i'} C_i \| x \mid K_n \|_n \\
& < & (1 + \frac{1}{\ell})|T_{i_\subl}^* x_{i_\subl, \ell}^*(x) |\\
& \leq & \norm x \norm.
\end{eqnarray*}
Also, for every $i \in I$ with $i \leq i'$ there exists $\ell(i) \in
\N$ such that $i_\sub{\ell(i)}=i$ and 
\[ \| x \mid K_i \|_i = | x_{i_\sub{\ell(i)}, \ell(i)}^* (x \mid
K_i) |= |(T_{i_\sub{\ell(i)}}^* x_{i_\sub{\ell(i)},
\ell(i)}^*)(x)| \]
Hence
\begin{eqnarray*}
\max_{ i \in I, i \leq i'} \| x \mid K_i \|_i & = & \max_{i \in I, i
\leq i'} | ( T_{i_\sub{\ell(i)}}^* x_{i_\sub{\ell(i)},
\ell(i)}^*)(x)| \\
& < & \max_{i \in I, i \leq i'}(1 + \frac{1}{\ell(i)})|(
T_{i_\sub{\ell(i)}} ^* x_{i_\sub{\ell(i)}, \ell(i)}^*)(x)| \\
& \leq & \norm x \norm.
\end{eqnarray*}
Therefore in case 1 we have that $| x |< \norm x \norm$.

\bigskip
\noindent
\underline{Case 2:} Assume that $n \in J$.

\bigskip
\noindent
Since $\| x \| = \| x \mid K_n \|_\infty$ we obtain 
\[ \| x \mid K_n \|_\infty \leq \frac{1}{1- \e_n} \max_{a \in A_n}
|a(x) |. \]
Hence
\begin{eqnarray*}
\sup_{j \in J} \max_{a \in A_j} |a(x) | & \leq & \sup_{j \in J} \| x \mid
K_j \|_\infty \\
& \leq & \| x \mid K_n \|_\infty \leq \frac{1}{1- \e_n} \max_{a \in
A_n} | a(x) | \\
& < & (1 + \delta_n) \max_{a \in A_n} |a(x)| \leq \norm x \norm.
\end{eqnarray*}
Also, if $\{ i \in I: i > n \} \not = \emptyset$ then we have that 
\begin{eqnarray*}
\sup_{i \in I, i >n} \| x \mid K_i \|_i & \leq & \sup_{i \in I, i>n}
C_i \| x \mid K_i \|_\infty \\
& \leq & \sup_{i \in I, i>n} C_i \| x \mid K_n \|_\infty \\
& \leq & (\sup_{i \in I, i>n}C_i) \frac{1}{1-\e_n} \max_{a \in A_n}
|a(x)| \\
& < & (1+ \delta_n) \max_{a \in A_n} |a(x)| \\
& \leq & \norm x \norm.
\end{eqnarray*}
Finally, the argument in the final part of Case 1, shows that 
\[ \sup_{i \in I, i \leq n} \| x \mid K_i \|_i < \norm x \norm. \]
Therefore, in Case 2 we have that $| x | < \norm x \norm$. Hence, we
have shown that $| x | < \norm x \norm$ for all $ x \in X \backslash
\{ 0 \}$.

\bigskip
\noindent
Set 
\[ A = \{ \pm (1 + \frac{1}{\ell})T_{i_\subl}^*x_{i_\subl , \subl}^* :
\ell \in \N \} \cup \{ \pm (1 + \delta_j)a: j \in J, a \in A_j \}. \]
We now show that 
\begin{equation} \label{E:next}
\next Ba(X, \norm \cdot \norm)^* \subseteq A. 
\end{equation}
This will finish the proof of the lemma, since $A$ is countable.

\bigskip
\noindent
By the definition of $\norm \cdot \norm$ we have that the set of the
extreme points of the dual ball, $\ext Ba(X, \norm \cdot \norm)^*$
satisfies
\[ \ext Ba(X, \norm \cdot \norm )^* \subseteq \wstarcl (A) \]
where for a set $B \subseteq X^*$, $\wstarcl (B)$ denotes the weak$^*$
closure of $B$. Let also $\wstarap (B)$ to denote the set of the
proper weak$^*$ accumulation points of $B$, i.e.
\[ \wstarap (B) = \wstarcl (B) \backslash B. \]
We shall show that if $x^* \in \wstarap (A)$ and there exists $x \in
X$ with $\norm x \norm =1$ such that $\norm x^* \norm = | x^*(x) |$,
then $\norm x^* \norm < 1$. Thus $x^* \not \in \ext Ba(X, \norm \cdot
\norm)^*$ which gives (\ref{E:next}).
Consider $x^* \in \wstarap (A)$ and $x \in X$ with $\norm x \norm =1$
such that $\norm x^* \norm = |x^*(x)|$. Then we have that 
\[ \norm x^* \norm = | x^* (x) | \leq | x^*| |x| < |x^* | \norm x
\norm = | x^* |. \]
Since $x^* \in \wstarap (A)$ we obtain that 
\[ x^* \in \wstarcl [ \{ \pm (1 + \frac{1}{\ell})T_{i_\subl}^*
x_{i_\subl , \subl}^*: \ell > \ell ' \} \cup \{ \pm (1+ \delta_j) a: J
\ni j > j',\; a \in A_j \}], \]
for all $\ell ' \in \N$ and for all $j' \in J$.
Since $\lim_{J \ni j \rightarrow \infty} \delta_j =0$ (if $J$ is
infinite) we have that $|x^*| \leq 1$. Thus $\norm x ^*\norm < 1$ which
finishes the proof. \hfill $\Box$

\section{The proof of the main result} \label{S:proof}

Lemma \ref{L:ip} is one of the basic ingredients for the proof of
Theorem \ref{T:main}. Before we present the proof of Theorem
\ref{T:main}, we give some more preliminary ingredients. We use the
following subsequence dichotomy for the $c_0$ basis, due to J. Elton:

\begin{Thm}[\cite{E1}] \label{T:c0dichotomy}
Every semi-normalized weakly null sequence  which does not have a
semi-boundedly complete subsequence, has a subsequence
equivalent to the unit vector basis of $c_0$.
\end{Thm}

\noindent
Recall that a sequence $(x_n)$ is called {\it semi-boundedly complete}
if for every sequence $(\lambda_n) \subset \R$ we have
\[ \sup_m \| \sum_{n=1}^m \lambda_n x_n \| < \infty \im \lambda_n
\rightarrow 0. \]
Our main result will follow from the following:

\begin{Thm} \label{T:premain}
If $(x_n)$ is a basic sequence in a Banach space $X$ with $\inf_n \|
x_n \| >0$, 
and $E$ is an isomorphically precisely norming set for $X$ such that 
\[ \sum_n | x^*(x_{n+1}-x_n) | < \infty, \; \forall x^* \in E, \]
then there exists a subsequence of $(x_n)$ which spans an
isomorphically polyhedral Banach space. 
\end{Thm}

\noindent
We postpone the proof of Theorem \ref{T:premain} for the moment. We
first give a proof of Theorem \ref{T:main} using the result of Theorem
\ref{T:premain}. First of all we present a proof of a well known
result (Lemma \ref{L:basic}) which generalizes the Bessaga-Pelczynski
Selection Principle \cite{BP} (see also pages 42 and 46 in \cite{D}).

\begin{Def} Let $(X, \| \cdot \|)$ be a Banach space and $Y$ be a linear subspace
(not necessarily closed) of $X^*$. 
\begin{itemize}
\item[(a)] $Y$ is a total linear space if for every $x \in X$ we have
that:
\begin{quote}
If $y(x)=0$ for all $y \in Y$ then $x=0$.
\end{quote}
\item[(b)] $Y$ is a 1-norming linear space if for every $x \in X$ we
have that 
\[ \| x \| = \sup_{y \in Y, \| y \|=1} |y(x)|. \]
\end{itemize}
Also, for a functional $x^* \in X^*$ we define the kernel of $x^*$ by:
\[ \ker x^* = \{ x \in X : x^*(x)= 0 \}. \]
\end{Def}

\noindent
The proof of the Bessaga-Pelczynski Selection Principle gives that if
$(x_n)$ is a normalized sequence in a Banach space and $Y$ is a
1-norming linear (not necessarily closed) space such that 0 is a
cluster point of $(x_n)$ for the $\sigma (X,Y)$ topology, then $(x_n)$
admits a basic subsequence. 

\begin{Lem} \label{L:basic}
Let $(X, \| \cdot \|)$ be a Banach space and $Y$ be a total linear
(not necessarily closed) subspace of $X^*$ Let $(x_n)$ be a sequence
in $X$ such that $\inf_n \| x_n \| >0$ and $y(x_n) \rightarrow 0$ for
all $y \in Y$. Then $(x_n)$ has a basic subsequence.
\end{Lem}

\noindent
{\bf Proof} Set $y_n = x_n / \| x_n \|$ for all $n \in \N$. Let
$x^{**} \in X^{**}$ be a weak$^*$ cluster point of $(y_n)$. Since
$\inf_n \| x_n \| >0$, we have that $y(y_n) \rightarrow 0$ for all $y
\in Y$. Thus $x^{**}
\mid Y =0$. Therefore we obtain that $x^{**} \in (X^{**} \backslash X)
\cup \{ 0 \}$. Define ${\cal N} = \ker x^{**}$  and 
\[ \norm x \norm = \sup \{ x^*(x) : x^* \in {\cal N}, \| x^* \|=1 \},
\; \forall x \in X. \]
Evidently $\norm \cdot \norm$ is an equivalent norm on $X$ and ${\cal
N}$ is a 1-norming space for $(X, \norm \cdot \norm)$. Also, 0 is a
cluster point of $(y_n / \norm y_n \norm)$ for the $\sigma (X, {\cal
N})$ topology. The proof of Bessaga-Pelczynski Selection Principle
gives that $(y_n / \norm y_n \norm)$ (and therefore $(x_n)$ itself)
admits a basic subsequence. \hfill $\Box$

\bigskip
\noindent
{\bf Proof of Theorem \ref{T:main}} 
Let $(x_n)$ be a sequence in a Banach space $X$ which does not
converge in norm, and let $E$ be an isomorphically precisely norming
set for $X$ such that (\ref{E:main}) holds. We define the (not
necessarily  closed) subspace $Y= \spn (E)$ of
$X^*$. Then $Y$ is a total linear space. Also, for every $y \in Y$ the
sequence $(y(x_n))_n$ is Cauchy. Define the linear function $ f : Y
\longrightarrow \R$ by
\[ f(y)= \lim_n y(x_n), \; \forall y \in Y. \]
Define a (not necessarily closed) 1-codimensional subspace ${\cal N}$
of $Y$ by
\[ {\cal N}= \{ y \in Y : f(y)=0 \}. \]
We separate two cases:

\bigskip
\noindent
\underline{Case 1:} Assume that ${\cal N}$ is total.

\bigskip
\noindent
Since $(x_n)$ is
not norm convergent we can assume by considering an appropriate
subsequence that $\inf_n \| x_n \| >0$. Since  $y(x_n) \rightarrow 0,
\;  \forall y \in {\cal N}$ there exists a subsequence of $(x_n)$
which is  basic.
The result follows from Theorem \ref{T:premain}.

\bigskip
\noindent
\underline{Case 2:} Assume that ${\cal N}$ is not total.

\bigskip
\noindent
Since ${\cal N}$ is a non-total 1-codimensional linear subspace of $Y$,
and $f$ is a linear function on $Y$, there exists $x \in X$ such that 
\[ y(x)=f(y), \; \forall y \in Y. \]
Thus we have that 
\[ y(x_n -x) \rightarrow 0, \; \forall y \in Y. \]
Since $(x_n)$ is not norm convergent and $Y$ is total, considering an
appropriate subsequence , we can assume that $\inf_n \| x_n -x \| >0$,
$(x_n-x)$ is basic, and  
\[ \sum_n | x^* [(x_{n+1}-x)-(x_n-x) ] | < \infty, \; \forall x^* \in
E. \]
Thus by Theorem \ref{T:premain} there exists  a subsequence $(x_{n_k})$
of $(x_n)$ such that $[(x_{n_k} -x)_k]$ is isomorphically polyhedral.
Thus the 1-dimensional extension $[(x_{n_k}-x)_k]+[x]$ is an i.p.
space, and therefore so is its subspace $[x_{n_k}]$. 

\bigskip 
\noindent
Conversely, consider a separable isomorphically polyhedral Banach
space $Y$. By Theorem \ref{T:charactpolyh} there exists a countable
isomorphically precisely norming set $E= \{ f_1, f_2, \ldots \}$ of
non zero functionals. Also, since $Y$ is separable, there exists a
sequence $(y_n)$ such that $[y_n] = Y$. By rearranging (and perhaps by
increasing) the set of $y_i$'s we can assume that 
\begin{eqnarray*}
& & [y_1] + \ker f_1 = X, \\
& & [y_1, y_2] + \cap_{j=1}^2 \ker f_j = X, \\
& & [y_1, y_2, y_3] + \cap_{j=1}^3 \ker f_j = X \mbox{ etc.}
\end{eqnarray*}
Set $x_1 = y_1$. We can write $y_2 = y_2' + y_2''$ where $y_2' \in
[y_1]$ and $y_2' \in \ker f_1$. Set $x_2 = y_2 ''$. Again, we can
write $y_3 = y_3' + y_3 ''$ where $y_3' \in [y_1, y_2]$ and $y_3 ''
\in \cap_{j=1}^2 \ker f_j$. Set $x_3= y_3 ''$ etc. Therefore, for
every positive integer $n$ we have that 
\begin{eqnarray*}
[x_1, \ldots, x_n]= [y_1, \ldots, y_n] \mbox{ and} \\
x_{n+1} \in \cap_{j=1}^n \ker f_j.
\end{eqnarray*}
Thus, $[x_n]=Y$ and (\ref{E:main}) holds. \hfill $\Box$

\bigskip
\noindent
We now present the 

\bigskip
\noindent
{\bf Proof of Theorem \ref{T:premain}} We can assume without loss of
generality that $X$ is separable (e.g. by considering $X=[x_n]$). For
every $x \in X$ we define 
\[ \| x \| = \sup_{e \in E} |e(x)|. \]
This defines an equivalent norm on $X$, and $E$ is a precisely norming
set for $(X, \| \cdot \|)$. Also, the weak$^*$ topology is metrizable
on $Ba(X^*)$, and let $d( \cdot , \cdot)$ denote the induced metric.
For $m \in \N$ we define (set $x_0=0$)
\[ K_m = \{ x^* \in Ba(X, \| \cdot \|)^*: \sum_{n=1}^\infty
|x^*(x_{n}-x_{n-1})| \leq m \}. \]
Then, $K_m$ is a weak$^*$ closed subset of $Ba(X^*)$ for every $m \in
\N$, $K_1 \subseteq K_2 \subseteq \cdots$, and $K:= \cup_{m=1}^\infty
K_m \supseteq E$. Define $f:K \longrightarrow \R$ by 
\[ f(k)= \lim_n k(x_n), \; \forall k \in K. \]
We separate the following cases:

\bigskip
\noindent
\underline{Case 1:} Assume that there exists $m \in \N$ such that the
restriction $f \mid K_m$ is not continuous ($K_n$ will always be equipped
with the weak$^*$ topology of $X^*$, for every $n \in \N$).

\bigskip
\noindent
{\bf Claim A:} For every $m' \geq m$ there exists a subsequence
$(x_n^{m'})_n$ of $(x_n)$ satisfying:
\begin{itemize}
\item $(x^m_n)_n$ is a subsequence of $(x_n)$.
\item $(x^{m'+1}_n)_n$ is a subsequence of $(x_n^{m'})_n$.
\item $[(x^{m'}_n \mid K_{m'})_n]$ is an i.p. Banach space (where
$[(x_n^{m'} \mid K_{m'})_n]$ denotes the completion of the normed
space $\spn (x_n^{m'} \mid K_{m'})_n$). 
\end{itemize}

\bigskip
\noindent
Indeed, for $m'=m$ we have that 
\[ \sup \{ \sum_n | x^*(x_n-x_{n-1})|: x^* \in Ba([(x_n \mid K_m)_n],
\| \cdot \|_{C(K_m)})^* \} \leq m \]
and $(x_n \mid K_m)_n$ is non-weakly convergent in $C(K_m)$. Thus by
Theorem \ref{T:bp} there exists a subsequence $(x^m_n)_n$ of $(x_n)$
such that $(x^m_n \mid K_m)_n$ is equivalent to the summing basis.
Thus $[(x^m_n \mid K_m)_n]$ is an i.p. Banach space. The proof of the
inductive step is a repetition of the same argument, since the
hypothesis ``$f \mid K_m$ is not continuous'' gives that ``$f \mid
K_{m'}$ is not continuous'' for every $m' \geq m$. The proof of Claim
A is complete.

\bigskip
\noindent
Let $(C_{m'})_{m' \geq m} \subset (1,2)$ be a sequence of numbers with
$C_{m'} \searrow 1$. Set $y_n=x^n_n$ for every $n \geq m$. Then
$(y_n)_{n \geq m}$ is a subsequence of $(x_n)$, and $(y_n)_{n \geq m'}$
is a subsequence of $(x^{m'}_n)_n$ for every $m' \geq m$. Thus $[(y_n
\mid K_{m'})_{n \geq m'}]$ is an i.p. Banach space for every $m' \geq m$.
Hence $[(y_n \mid K_{m'})_{n \geq m}]$ is an i.p. Banach space for every $m'
\geq m$.
Thus, by Theorem \ref{T:dfh} and Theorem \ref{T:charactpolyh}, for every
$m' \geq m$ there exists a $\sqrt{C_{m'}}$ equivalent norm $\| \cdot
\|_{m'}$ on $[(y_n \mid K_{m'})_{n \in \N}]$ such that 
\[ \next Ba ([(y_n \mid K_{m'})_{n \geq m}], \| \cdot \|_{m'})^*
\mbox{ is countable}. \]
Also, by multiplying $\| \cdot \|_{m'}$ by an appropriate constant, we
can assume that 
\[ \| z \|_{C(K_{m'})} \leq \| z \|_{m'} \leq C_{m'} \| z
\|_{C(K_{m'})}, \; \forall z \in [(y_n \mid K_{m'})_{n \geq m}]. \]
Thus the assumptions of Lemma \ref{L:ip} are satisfied for the Banach
space $[y_n]$, $I=\{ m, m+1, \ldots \}$, and $J = \emptyset$. Hence
$[y_n]$ is an i.p. space.

\bigskip
\noindent
\underline{Case 2:} Assume that $f \mid K_m$ is continuous for every
$m \in \N$.

\bigskip
\noindent
We separate two cases:

\bigskip
\noindent
\underline{Subcase 2.1:} Assume that there exists a subsequence $(y_n)$
of $(x_n)$ and there exists $m \in \N$ such that 
\begin{eqnarray*}
\inf_{n \not = n'} \| (y_n - y_{n'}) \mid K_m \|_{C(K_m)} >0.
\end{eqnarray*}

\noindent
Hence 
\begin{eqnarray*}
\inf_n \| (y_n-f) \mid K_m \|_{C(K_m)} > 0
\end{eqnarray*}
and therefore for every $m' \geq m$ we have that 
\begin{eqnarray*}
\inf_n \| (y_n -f) \mid K_{m'} \|_{C(K_{m'})} >0.
\end{eqnarray*}
Thus for every $m' \geq m$, $((y_n-f) \mid K_{m'})_n$ is a weakly null
semi-normalized sequence (by the definition of $K_{m'}$, note that $\|
y_n \mid K_{m'} \|_{C(K_{m'})} \leq m', \; \forall n \in \N$).

\bigskip
\noindent
{\bf Claim B:} For every subsequence $(z_n)$ of $(y_n)$ and for every
$m'=m, m+1, \ldots$ we have that $((z_n -f) \mid K_{m'})_n$ is not
semi-boundedly complete.

\bigskip
\noindent
Indeed, for every $n \in \N$ we have that 
\[ \| [(z_1-f)-(z_2-f)+ \cdots +(-1)^{n+1}(z_n -f)] \mid K_{m'}
\|_{C(K_{m'})} \]
\[\leq \| [z_1-z_2+ \cdots +(-1)^{n+1}z_n] \mid K_{m'} \|_{C(K_{m'})}
+ \| f \mid K_{m'}\|_{C(K_{m'})}.\]
There exists $k \in K_{m'}$ such that 
\begin{eqnarray*}
& & \| [z_1-z_2+ \cdots +(-1)^{n+1}z_n ] \mid K_{m'} \|_{C(K_{m'})}\\ 
& & \hskip 2in = | (z_1-z_2+ \cdots + (-1)^{n+1}z_n)(k)| \\
& & \hskip 2in \leq   |(z_1-z_2)(k)| + |(z_3-z_4)(k)| + \cdots  + m'\\
& & \hskip 2in \leq  \sum_i |k(x_i-x_{i-1})| +m' \\
& & \hskip 2in \leq  2 m'.
\end{eqnarray*}
Thus 
\begin{eqnarray*}
& & \sup_n \| [(z_1-f)-(z_2-f)+ \cdots + (-1)^{n+1}(z_n-f)] \mid K_{m'}
\|_{C(K_{m'})} \\
& & \hskip 2in \leq 2m' + \|f \mid K_{m'} \|_{C(K_{m'})}.
\end{eqnarray*}
Therefore, the sequence $((z_n-f) \mid K_{m'})_n$ is not
semi-boundedly complete since the sequence $((-1)^{n+1})_n$ does not
converge to zero. Claim B is proved.

\bigskip
\noindent
{\bf Claim C:} For every $m' \geq m$ there exists a subsequence
$(y^{m'}_n)_n$ of $(y_n)$ satisfying
\begin{itemize}
\item $(y^m_n)$ is a subsequence of $(y_n)$.
\item $(y^{m'+1}_n)_n$ is a subsequence of $(y^{m'}_n)_n$.
\item $([(y^{m'}_n \mid K_{m'})_n], \| \cdot \|_{C(K_{m'})})$ is an
i.p. Banach space.
\end{itemize}

\bigskip
\noindent
Indeed, for $m'=m$, $((y_n-f) \mid K_m)_n$ is a weakly null semi-normalized
sequence which does not have any semi-boundedly complete subsequence
(by Claim B). By Theorem \ref{T:c0dichotomy} there exists a
subsequence $(y^m_n)_n$ of $(y_n)$ such that $((y^m_n -f) \mid K_m)_n$
is equivalent to the unit vector basis of $c_0$. Thus $([((y^m_n-f)
\mid K_m)_n], \| \cdot \|_{C(K_m)})$ is an i.p. Banach space. Hence
$[((y^m_n-f) \mid K_m)_n]+[f \mid K_m]$ is an i.p. Banach space, and
therefore so is its subspace $[(y^m_n \mid K_m)_n]$. The proof of the
inductive step is a repetition of the same argument. The proof of Claim C
is  complete and the proof
of Subcase 2.1 finishes identically as in Case 1.

\bigskip
\noindent
\underline{Subcase 2.2:} Assume that for every subsequence $(y_n)$ of
$(x_n)$, and for every $m \in \N$ we have that 
\[ \inf_{n \not = n'} \|(y_n -y_{n'})|K_m \|_{C(K_m)}=0. \]

\noindent
Let $(\e _n') \subset (0,1)$ be a sequence of numbers such that $n \e_n'
\searrow 0$.

\bigskip
\noindent
{\bf Claim D:} There exists a subsequence $(y_n)$ of $(x_n)$ such that

\[ \sum_{n=m}^\infty \| (y_n -f) \mid K_m \|_{C(K_m)} \leq \e_m', \;
\forall m \in \N. \]

\noindent
In order to prove Claim D, we first prove

\bigskip
\noindent
{\bf Subclaim Da:} For every subsequence $(y_n)$ of $(x_n)$, for every
$m \in \N$, and for every $\e >0$ there exists a subsequence $(z_n)$
of $(y_n)$ such that 
\[ \| (z_1-z_n) \mid K_m \|_{C(K_m)} < \e, \; \forall n \in \N. \]

\bigskip
\noindent
Indeed, assume that Subclaim Da is false. Thus, if we set
\[ I_1= \{ n \in \N: \|(y_1 -y_n) \mid K_m \|_{C(K_m)} < \e \}, \]
then $I_1$ is finite. Set $i_1 = \max I_1 +1$. Also, the set 
\[ I_2 = \{ n >i_1 : \| (y_{i_1} -y_n ) \mid K_m \|_{C(K_m)} < \e \}\]
is finite. Set $i_2 = \max I_2 +1$. We continue similarly. Then the
subsequence $(y_{i_n})$ of $(y_n)$ satisfies 
\[ \inf_{n \not = n'} \| (y_{i_n}-y_{i_{n'}}) | K_m \|_{C(K_m)} > \e \]
which is a contradiction. Subclaim Da is proved.

\bigskip
\noindent
{\bf Subclaim Db:} For every subsequence $(y_n)$ of $(x_n)$ and for
every $m \in \N$ there exists a subsequence $(z_n)$ of $(y_n)$ such
that 
\[ \sum _{n=1}^\infty \| (z_n -f) \mid K_m \|_{C(K_m)} \leq \e_m'. \]

\noindent
Indeed, using Subclaim Da, choose a subsequence $(y^1_n)$ of $(y_n)$
such that 
\[ \| (y^1_1 -y^1_n) \mid K_m \|_{C(K_m)} \leq \frac{\e _m'}{2^2}, \;
\forall n \in \N. \]
Then, we choose a subsequence $(y^2_n)$ of $(y^1_n)$ such that 
\[ \|(y^2_1 -y^2_n) \mid K_m \|_{C(K_m)} \leq \frac{\e_m'}{2^3}, \;
\forall n \in \N. \] 
We continue similarly to define $(y^k_n)_n$ for every $k \in \N$
satisfying:
\[ \| (y^k_1 - y^k_n) \mid K_m \|_{C(K_m)} \leq \frac{\e_m'}{2^{k+1}},
\; \forall n \in \N. \]
Thus
\[ \| (y^k_n -y^k_{n'}) \mid K_m \|_{C(K_m)} \leq \frac{\e_m'}{2^k}, \;
\forall n, n' \in \N. \]
Set $z_n=y^n_n, \; \forall n \in \N$. Then $(z_n)$ is a subsequence of
$(y_n)$ and for every $k \in \N$ $(z_n)_{n \geq k}$ is a subsequence
of $(y^k_n)_n$. Hence
\[ \| (z_n - z_{n'}) \mid K_m \|_{C(K_m)} \leq \frac{\e_m'}{2^k}, \;
\forall k \in \N, \; \forall n,n' \geq k. \]
By taking $n=k$ and $n' \rightarrow \infty$ we obtain that 
\[ \| (z_n-f) \mid K_m \|_{C(K_m)} \leq \frac{\e_m'}{2^n}, \; \forall n
\in \N, \]
which proves Subclaim Db.

\bigskip
\noindent
Now the proof of Claim D consists of another diagonal argument using
Subclaim Db: Choose a subsequence $(x_n^1)$ of $(x_n)$ such that 
\[ \sum_{n=1}^\infty \| (x^1_n -f) \mid K_1 \|_{C(K_1)} \leq \e_1'. \]
Then, choose a subsequence $(x^2_n)$ of $(x^1_n)$ such that 
\[ \sum_{n=1}^\infty \| (x^2_n -f) \mid K_2 \|_{C(K_2)} \leq \e_2' . \]
We continue similarly to choose the sequences $(x^k_n)_n$ for every $k
\in \N$,   and then we define the subsequence $(y_n)$ of $(x_n)$ by
$y_n=x^n_n, \; \forall n \in \N$. Since for every $m \in \N$ we have
that $(y_n)_{n \geq m}$ is a subsequence of $(x^m_n)_n$, the
conclusion of Claim D follows. 

\bigskip
\noindent
We shall prove that $[y_n]$ is an i.p. space.

\bigskip
\noindent
Since $(y_n)$ is a basic sequence with $\inf_n \| y_n \|>0$, there
exists  $M >0$ such that if $\sum_{n=1}^\infty \lambda_n y_n \in
[y_n]$  with $\| \sum_{n=1}^\infty \lambda_n y_n \| \leq 1$ then
\begin{itemize}
\item $\| \sum_{n=m}^\infty \lambda_n y_n \| \leq M, \; \forall m \in \N$, and 
\item $| \lambda_n | \leq M,\; \forall n \in \N$.
\end{itemize}

\noindent
{\bf Claim E:} Assume that for some $m_0 \in \N$ and for some $x \in
K_{m_0}$ we have that $f(x) \not =0$. Let $y = \sum_{n=1}^\infty
\lambda_n y_n \in [y_n]$ with $\| y_n \|=1$. Then $\sum_{n=1}^\infty
\lambda_n$ converges and 
\[ | \sum_{n=m}^\infty \lambda_n | \leq \frac{2M}{|f(x)|}, \;
\forall m \geq m_0. \]

\bigskip
\noindent
Indeed, for every $m_0 \leq m \leq m'$ we have that 
\begin{eqnarray*}
| \sum_{n=m}^{m'} \lambda_n y_n(x) - \sum_{n=m}^{m'} \lambda_n f(x)| &
\leq & \sum_{n=m}^{m'} | \lambda_n | | y_n(x) - f(x) | \\
& \leq & M \sum_{n=m}^{m'} \| (y_n -f) \mid K_{m_0} \|_{C(K_{m_0})}. 
\end{eqnarray*}
Thus $\sum_{n=1}^\infty \lambda_n$ converges since $\sum_{n=1}^\infty
\lambda_n y_n(x)$ converges and $\sum_{n=1}^\infty \| (y_n -f) \mid
K_{m_0} \|_{C(K_{m_0})}$ converges. Also, by taking $m' \rightarrow
\infty$ we obtain that 
\begin{equation} \label{E:closetails}
| \sum_{n=m}^\infty \lambda_n y_n(x) - \sum_{n=m}^\infty \lambda_n
f(x) | \leq M \e_m' \leq M, \; \forall m \geq m_0.
\end{equation}
We also have that 
\begin{equation} \label{E:boundedtails}
\| \sum_{n=m}^\infty \lambda_n y_n \| \leq M. 
\end{equation}
Combining (\ref{E:closetails}) and (\ref{E:boundedtails}) we obtain
that
\[ | \sum_{n=m}^\infty \lambda_n | = \frac{1}{|f(x)|} |
\sum_{n=m}^\infty \lambda_n f(x)| \leq \frac{2M}{|f(x)|}. \]
The proof of Claim E is complete.

\bigskip
\noindent
Now, if $f \mid K \not = 0$ then choose $m_0 \in \N$ and $x \in
K_{m_0}$ such that $f(x) \not = 0$. Otherwise set $m_0=1$.

\bigskip
\noindent
The rest of the proof is similar to the proof of Theorem 1 in
\cite{F4}.

\bigskip
\noindent
For every $m \geq m_0$ the set of functions $\{ y_1 \mid K_m, \ldots ,
y_m  \mid K_m , f \mid K_m \}$ is uniformly continuous on $(K_m,$
 weak$^*$  topology). Therefore, there exists $\delta_m >0$ such that 
\begin{equation} \label{E:unifcont}
k_1, k_2 \in K_m \; \& \; d(k_1, k_2) < \delta_m \im \left\{
\begin{array}{l} 
| y_n(k_1) - y_n(k_2) | < \e_m', \; \forall n=1, \ldots , m \\
\mbox{and } |f(k_1)-f(k_2)|< \e_m'. \end{array} \right.
\end{equation}
For every $m \geq m_0$ since $(K_m, d( \cdot , \cdot ))$ is a compact
metric space, we can choose a finite set $A_m \subset K_m$ such that 
\begin{equation} \label{E:net}
\forall k \in K_m \; \exists a \in A_m \; d(k,a) < \delta_m . 
\end{equation}

\noindent
{\bf Claim F:} For $y \in [y_n]$ and $m \geq m_0$ such that $1=\| y\| =
\| y \mid K_m \|_{C(K_m)}$, there exists $a \in A_m$ such that 
\[ | y(a)| \geq  \left\{ \begin{array}{ll}
1-\frac{5mM}{|f(x)|} \e_m' & \mbox{ if } f \not = 0\\
1-3mM \e_m' & \mbox{ if } f = 0.
\end{array} \right. \]

\bigskip
\noindent
Indeed, let $y= \sum_n \lambda_n y_n$ with $1= \| y \|= \| y \mid K_m
\|_{C(K_m)}$ and let $k \in K_m$ with $|y(k)|=1$. 
By (\ref{E:net}), there exists $a \in A_m$ with
$d(k,a)< \delta_m$. Thus by (\ref{E:unifcont}) we have that 
\begin{equation} \label{E:unifcontinuity}
|f(a)-f(k)|< \e_m' \mbox{ and } |y_n(a)-y_n(k)|< \e_m' \mbox{ for }n=1,
\ldots,m.
\end{equation}
The following computation is valid whether $f \not =0$ or $f=0$. If
$f=0$ then we set $(\sum_{n=m+1}^\infty \lambda_n f) \mid K_m =0$.
\begin{eqnarray*}
| y(a) -y(k) |& = & | \sum_{n=1}^\infty \lambda_n y_n(a) -
\sum_{n=1}^\infty \lambda_n y_n(k) | \\
& \leq & \sum_{n=1}^m | \lambda_n | | y_n (a) -y_n(k) | + |
\sum_{n=m+1}^\infty \lambda_n y_n(a) - \sum_{n=m+1}^\infty \lambda_n
y_n (k) | \\
& \leq & m M \e_m' + | \sum_{n=m+1}^\infty \lambda_n y_n (a)
-\sum_{n=m+1}^\infty \lambda_n f(a) + \sum_{n=m+1}^\infty \lambda_n
f(a) \\
& & - \sum_{n=m+1}^\infty \lambda_n f(k) + \sum_{n=m+1}^\infty
\lambda_n f(k) - \sum_{n=m+1}^\infty \lambda_n y_n(k) | \\
& \leq & mM \e_m' + \sum_{n=m+1}^ \infty | \lambda_n ||y_n(a) - f(a) |
\\
& & + | \sum_{n=m+1}^\infty \lambda_n | |f(a)-f(k) | +
\sum_{n=m+1}^\infty | \lambda_n | | f(k) -y_n(k) |.
\end{eqnarray*}
Thus if $f \not = 0$, by (\ref{E:unifcontinuity}), and Claims D and E
we obtain that 
\begin{eqnarray*}
| y(a) -y(k) | & \leq & mM\e_m' + M\e_m' + \frac{2M}{|f(x)|} \e_m' + M
\e_m' \\
& \leq & \frac{5mM}{|f(x)|}\e_m'.
\end{eqnarray*}
If $f=0$, by (\ref{E:unifcontinuity}) and Claim D we obtain that
\[ | y(a)-y(k) | \leq mM \e_m' + M \e_m' + M \e_m' \leq 3mM \e_m'. \] 
Since  $|y(k)|=1$, the proof of Claim F is complete.

\bigskip
\noindent
Choose $m_0' \geq m_0$ such that 
\[ \begin{array}{ll}
1-\frac{5nM}{|f(x)|}\e_n' \in (0,1) & \mbox{ if } f \not =0, \mbox{ and}\\
1-3nM \e_n' \in (0,1) & \mbox{ if } f=0
\end{array} \]
for all $n \geq m_0'$. Then choose a sequence $(\e_n)_{n \geq m_0'}
\subset (0,1)$ with $\lim_n \e_n =0$, such that 
\[ \begin{array}{ll}
\frac{1}{1-\e_n}(1- \frac{5nM}{|f(x)|} \e_n')>1 & \mbox{ if }f \not
=0, \mbox{ and}\\
\frac{1}{1-\e_n}(1-3nM\e_n')>1 & \mbox{ if }f=0
\end{array} \]
for all $n \geq m_0'$. 
Claim F shows that the assumptions of Lemma \ref{L:ip} are satisfied
for the Banach space $[y_n]$, $I= \emptyset$, $J= \{ m_0', m_0'+1,
\ldots \}$ and $(\e_n)_{n \geq m_0'}$. Thus $[y_n]$
is an i.p. space which finishes the proof of Theorem \ref{T:premain}.
\hfill $\Box$

\bigskip
\noindent
Using Theorem 1 of \cite{F4} we can give an easy proof of the
following weaker result than Theorem \ref{T:main}.

\bigskip
\noindent
{\bf Remark \ref{R:fonf}}
Under the same hypotheses of Theorem \ref{T:main} there exist a
sequence $(\e _n) \in \{ \pm 1 \} ^{\N}$ and an increasing sequence of
positive integers $(\ell_k)$ such that $[( \sum_{i=1}^{\ell_k} \e_i
(x_i - x_{i-1}))_k]$ is an i.p. space.

\bigskip
\noindent
Indeed, the proof of Theorem 1 in \cite{F4} shows the following:
\begin{quote}
Let $(X, \| \cdot \|)$ be a Banach space, $K_1 \subset K_2 \subset
\cdots$ be subsets of $Ba (X^*)$ and let $(w_n)$ be a sequence in
$X$. If $(w_n)$ is basic, $\inf_n \| w_n \| >0$, $\sum_n \| w_n \mid
K_n \|< \infty$ and $\cup_n K_n$ is an isomorphically precisely
norming set, then $[w_n]$ is an i.p. Banach space.
\end{quote}

\noindent
Now, the proof of the assertion of the remark can be sketched as
follows: If there is no subsequence of $(x_n)$ equivalent to the
summing basis, then there exists a sequence $(\e_n) \in \{ \pm 1
\}^{\N}$ such that 
\[ (\sum_{i=1}^n \e_i (x_i - x_{i-1}))_n \mbox{ is not bounded.} \]
Therefore there exists an increasing sequence $(n_k)$ of integers such
that 
\[ \| \sum_{i=1}^{n_k} \e_i (x_i - x_{i-1}) \| \geq 2^k k, \; \forall
k \in \N. \]
Set $z_k= \sum_{i=1}^{n_k} \e_i (x_i -x_{i-1})$ for every $k \in
\N$. Since $(z_k)$ does not converge in norm, and $(y(z_k))$ is Cauchy
for every $y \in \spn E$, we obtain (as in the proof of Theorem
\ref{T:main}) that there exists $z \in X$ ($z$
can also be zero) and an increasing sequence
$(m_k)$ of integers such that $(z_{m_k}-z)$ is a basic sequence. Set 
\[ K_m = \{ x^* \in Ba(X^*): \sum_{n=1}^\infty | x^*(x_n -x_{n-1})|
\leq m \}, \forall m \in \N \]
(where $x_0=0$). We easily see that 
\[ \sum_k \| \frac{z_{m_k}-z}{\| z_{m_k}-z \|} \mid K_k \| <
\infty. \]
Thus, by the above mentioned Theorem 1 of \cite{F4} we obtain that 
\[ [( \sum_{i=1}^{n_{m_k}} \e_i (x_i - x_{i-1}))_k] \mbox{ is an
i.p. space.} \]

\section{Applications} \label{S:applications}

\noindent
As a first application we strengthen a corollary of Theorem
\ref{T:elton} which was also proved in a different way by  R. Haydon,
E. Odell and H. Rosenthal \cite{HOR}. First we need some definitions.
Let $K$ be a compact metric space. $B_1(K)$ denotes the class of
bounded  Baire-1 functions on $K$, i.e. the pointwise limits of the
uniformly bounded sequences of continuous functions on $K$. $DSC(K)$
denotes the space of bounded Differences of Semi-Continuous functions
on $K$, i.e. 
\begin{eqnarray*}
DSC(K) & = & \{ f : K \longrightarrow \R \mid \mbox{ there exists a uniformly
bounded sequence} \\
& & (f_n)_{n=1}^\infty \subset C(K) \mbox{ such that } \lim_n f_n
(k)=f(k)  \mbox{ and }\\
& & \sum_{n=1}^\infty | f_{n+1}(k) - f_n(k)| < \infty \mbox{ for all
}k \in K \}. 
\end{eqnarray*}
Let $f$ be a non-continuous function on $B_1(K)$ and ${\cal C}$ be a
non-empty class of Banach spaces. Using terminology which was
introduced by R. Haydon, E. Odell and H. Rosenthal \cite{HOR}, we say
that {\it $f$ governs ${\cal C}$} if for every uniformly bounded
sequence $(f_n)$ of continuous functions on $K$ which converges
pointwise to $f$ on $K$, there exists $X \in {\cal C}$ which embeds
isomorphically in the closed linear span $[f_n]$ of $(f_n)$ equipped
with the supremum norm. We say that {\it $f$ strictly governs ${\cal
C}$} if for every uniformly bounded sequence $(f_n)$ of continuous
functions on $K$ which converges pointwise to $f$ on $K$ there exists
a convex block sequence $(g_n)$ of $(f_n)$ such that the closed linear
span $[g_n]$ of $(g_n)$ is isomorphic to some $X \in {\cal C}$. A
corollary of Theorem \ref{T:elton} which was proved in a different way
by R. Haydon, E. Odell and H. Rosenthal can be stated as follows:

\bigskip
\noindent
{\bf Theorem \ref{T:hordsc} [\cite{E2}, \cite{HOR}]}
Let $f \in DSC(K) \backslash C(K)$ be given, where $K$ is a compact
metric  space. Then $f$ governs $\{ c_0 \}$.

\bigskip
\noindent
A generalization of this result is the following:

\bigskip
\noindent
{\bf Theorem \ref{T:dsc}}
Let $f \in DSC(K) \backslash C(K)$ be given, where  $K$ is a compact
metric  space. Then $f$ strictly governs the class of (separable)
polyhedral  Banach spaces.

\bigskip
\noindent
For deducing Theorem \ref{T:dsc} from Theorem \ref{T:main} we need
the next well known remark. We first fix some terminology: If $A$ is a
subset of a Banach space $X$ then $\tilde{A}$ denotes the weak$^*$
closure of $A$ in $X^{**}$. Also if $A, B$ are non-empty subsets of
$(X, \| \cdot \|)$ then the minimum distance between $A$ and $B$ is
defined by:
\[{\rm md}\, (A,B)= \inf \{ \| a-b \| : a \in A, b \in B \}. \]

\begin{Rmk}\label{R:convexdistance}
If $A$, $B$ are convex subsets of a Banach space, then ${\rm md} \,
(A, B)= {\rm md}\, ( \tilde{A}, \tilde{B})$.
\end{Rmk}

\noindent
Thus, if $f \in DSC(K) \backslash C(K)$ and a bounded sequence $(f_n)$
of continuous functions which converges pointwise to $f$ on $K$, are
given, then by Remark \ref{R:convexdistance} there exists a convex
block sequence $(g_n)$ of $(f_n)$ such that 
\[ \sum_{n=1}^\infty | g_{n+1}(k) - g_n(k) | < \infty, \; \forall k
\in K. \]
Since $f \not \in C(K)$, we can also assume (by considering an
appropriate subsequence) that $(g_n)$ is a semi-normalized basic sequence.
Thus Theorem \ref{T:main} gives that some subsequence of $(g_n)$ spans
an i.p. Banach space, which proves Theorem \ref{T:dsc}.

\bigskip
\noindent
As a second application we obtain an Orlicz-Pettis type result:

\bigskip
\noindent
{\bf Theorem \ref{T:orliczpettis}}
Let $(y_n)$ be a sequence in a Banach space $X$ and let $E$ be an
isomorphically precisely norming set for $X$. If $c_0$ does not embed
isomorphically in the closed linear span $[y_n]$ of $(y_n)$ and 
\[ \sum_n | x^*(y_n) | < \infty , \; \forall x^* \in E, \]
then $\sum_n y_n$ converges unconditionally.

\bigskip
\noindent
{\bf Proof} For $(\eta_i) \subset \{ \pm 1 \} ^{\N} $ define the sequence
$(x_n)$ by 
\[ x_n = \sum_{i=1}^n \eta_i y_i, \; \forall n \in \N. \]
We have that the sequence $(x_n)$ satisfies (\ref{E:main}). Since
$c_0$ does not embed isomorphically in $[y_n]=[x_n]$, we have that the
conclusion of Theorem \ref{T:main} fails. Thus the sequence $(x_n)$
converges in norm. Hence $\sum_n y_n$ converges unconditionally. \hfill $\Box$
 
\bigskip
\noindent
As a final application of Theorem \ref{T:main} we prove the following
immediate corollary which has been proved previously by V. Fonf
\cite{F4}. 

\begin{Cor} 
Let $X$ be a Banach space which does not contain an
isomorph of $c_0$. Let $A$ be a subset of $X$, and let $B$ be an
isomorphically precisely norming subset of $X^*$. If for every $b \in
B$ the set $\{ b(a): a \in A \}$ is bounded, then $A$ is bounded.
\end{Cor}

\noindent
{\bf Proof} If $A$ is not bounded, we can find a sequence $(a_n)
\subset A$ such that $\| a_n \| > 2^n$ for all $n \in \N$. Set 
\[ \alpha_n = \sum_{i=1}^n \frac{a_i}{\| a_i \|}, \; \forall n \in
\N. \]
Thus 
\[ \sum | b( \alpha_{n+1}- \alpha_n )|< \infty, \; \forall b \in B. \]
Since $X$ does not contain an isomorph of $c_0$, by Theorem
\ref{T:main} we obtain that $(\alpha_n)$ converges in norm, which is a
contradiction. \hfill $\Box$

\noindent
G. Androulakis, Math. Sci. Bldg., University of Missouri-Columbia, 
Columbia, MO 65211-0001\\
e-mail: giorgis@math.missouri.edu 

\end{document}